\documentclass[11pt]{article}
\usepackage{amsfonts,latexsym}

\title{Non commutative two dimensional modular symbol}
\author{Ivan Emilov Horozov}
\date{September 16, 2006}

\newcommand{\beq}{\begin{equation}}
\newcommand{\eeq}{\end{equation}}
\newcommand{\beqa}{\begin{eqnarray}}
\newcommand{\eeqa}{\end{eqnarray}}
\newcommand{\beaa}{\begin{eqnarray*}}
\newcommand{\ben}{\begin{eqnarray*}}
\newcommand{\eaa}{\end{eqnarray*}}
\newcommand{\een}{\end{eqnarray*}}

\newcommand{\text}{\textrm}
\newcommand \nc {\newcommand}
\nc \proof {\noindent {\em{Proof.\/ }}} 
\nc \qed {$\Box$\hfill}
\newtheorem{theorem}{Theorem}[section]
\newtheorem{lemma}[theorem]{Lemma}
\newtheorem{proposition}[theorem]{Proposition}
\newtheorem{corollary}[theorem]{Corollary}
\newtheorem{definition}[theorem]{Definition}
\newtheorem{example}[theorem]{Example}
\newtheorem{remark}[theorem]{Remark}
\newtheorem{conjecture}[theorem]{Conjecture}
\newtheorem{question}[theorem]{Question}
\nc \bth[1] {\begin{theorem}\label{t#1} }
\nc \ble[1] {\begin{lemma}\label{l#1} }
\nc \bpr[1] {\begin{proposition}\label{p#1} }
\nc \bco[1] {\begin{corollary}\label{c#1} }
\nc \bde[1] {\begin{definition}\label{d#1}\rm }
\nc \bex[1] {\begin{example}\label{e#1}\rm }
\nc \bre[1] {\begin{remark}\label{r#1}\rm }
\nc \bcon[1] {\begin{conjecture}\label{con#1}\rm }
\nc \bque[1] {\begin{question}\label{que#1}\rm }
\nc {\eth} { \end{theorem} }
\nc {\ele} { \end{lemma} }
\nc {\epr}{ \end{proposition} }
\nc {\eco} { \end{corollary} }
\nc {\ede} {\end{definition} }
\nc {\eex} { \end{example} }
\nc {\ere} {\end{remark} } 
\nc {\econ} { \end{conjecture} }
\nc {\eque} {\end{question} }
\nc \eqref[1] {{\rm{(\ref{#1})}}} 
\nc \thref[1]{Theorem \ref{t#1}}
\nc \leref[1]{Lemma \ref{l#1}} 
\nc \prref[1]{Proposition
\ref{p#1}} \nc \coref[1]{Corollary \ref{c#1}} 
\nc \deref[1]{Definition \ref{d#1}} 
\nc \exref[1]{Example \ref{e#1}}
\nc \reref[1]{Remark \ref{r#1}} 
\nc \conref[1]{Conjecture\ref{con#1}}

\def \R {{\mathcal R}}

\def \K {{\mathcal K}}
\def\a{\alpha}
\def\b{\beta}

\def\ga{\gamma}
\def\Ga{\Gamma}

\def\si{\sigma}
\def\Si{\Sigma}

\def \Z {{\mathbb Z}}
\def \Q {{\mathbb Q}}
\def \R {{\mathbb R}}
\def \C {{\mathbb C}}

\nc \Wr {Wr} \nc \GRN { \Gr^{(N)} }
\nc \GRA[1] { \Gr_A^{(#1)} }   
\nc \GRAN { \GRA{N} } \nc \GrA[1] { \Gr_A(#1) }\nc \GrAa {
\GrA{\alpha} }
\nc \GRB[1] { \Gr_B^{(#1)} }   
\nc \GRBN { \GRB{N} } \nc \GrB[1] { \Gr_B(#1) } \nc \GrBb {
\GrB{\beta} }
\nc \GRMB[1] { \Gr_{MB}^{(#1)} }   
\nc \GRMBN { \GRMB{N} } \nc \GrMB[1] { \Gr_{MB}(#1) } \nc \GrMBb {
\GrMB{\beta} }

\begin{document}

\title{{\LARGE\bf{Non-commutative Two Dimensional Modular Symbol}}}

\author{
I. ~Horozov
\thanks{E-mail: ihorozov@brandeis.edu}
\\ \hfill\\ \normalsize \textit{Department of Mathematics,}\\
\normalsize \textit{ Brandeis University, 415 South St.,}\\
\normalsize \textit {MS 050, Waltham, MA 02454 }  \\ 
}
\date{}
\maketitle
\begin{abstract}
We define iteration over a two dimensional manifold as analog of
iteration over a path defined by Chen. We give several
applications. Some of them include constructions of non-abelian
modular symbol for $SL(3,\Z)$ and for $SL_{2/K}$, where $K$ is a real quadratic field.
We construct the latter as generating series of iterated integrals over a two-dimensional manifold of Hilbert modular forms. We give a motivic interpretation of this non-abelian
modular symbol. Other applications give
motivic interpretations of certain 
iterated completed Dedekind zeta functions.
\end{abstract}
\tableofcontents
\setcounter{section}{-1}
\section{Introduction}

In this paper we define iteration over a two dimensional manifold as analog of iteration over a path defined by Chen \cite{Ch}.We call this process iteration over a membrane. The main application of iteration over a two-dimensional membrane is the construction of non-abelian modular symbol for Hilbert modular forms associated to real quadratic fields. Another application is a construction of non-abelian modular symbol for $SL(3,\Z)$.  These modular symbols are generating series of iterated integrals over a two dimensional membrane. Modular symbol was first defined by Manin in \cite{M1}. Recently he defined a non-abelian modular symbol related to finite index subgroups of $SL(2,\Z)$ (see \cite{M2}). Manin used iteration over a path to define a non-abelian modular symbol for $SL_{2/\Q}$.

With the definition of iteration over a membrane we are able to give a construction of non-abelian modular symbol for Hilbert modular forms associated to real quadratic fields, previously not considered. We give a motivic interpretation of iterated integrals of Hilbert modular forms of real quadratic field over two dimensional membrane. In the last application we iterate certain Mellin transform of theta functions to define iterated Dedekind zeta function for some number fields $K$. For $K=\Q$ we obtain  multiple completed zeta functions. That is, we iterate a theta function. The theta function is such that if we take its Mellin transform we obtain Riemann zeta function together with the gamma factor. We give motivic interpretation of the  multiple completed zeta function when the variables take values in the positive even integers. For $K$ a real quadratic field of class number one we give motivic interpretation of the  multiple completed Dedekind zeta function when the var!
 iables take value in the positive even integers. For $K$ equal to an imaginary quadratic field of class number one or CM field of degree $4$ over $\Q$ of class number one we define and give motivic interpretation of multiple completed  Dedekind zeta function with variables taking value in the positive integers.

In the first section we define an iterated integral over a rectangle $A$. In certain cases we give a generalization to higher dimensions. We show that for given $n$ $2$-forms defined over a rectangle, the set of iterated integrals with $\Z$ coefficients have the structure of a Hopf algebra $R_A$. In order to show that the set of iterated integrals forms a Hopf algebra we assume that none of the iterated integrals over the rectangle vanishes. In case that some of the iterated integrals over the rectangle vanish then the $\Z$-module generated by the iterated integrals has the structure of a ring. And we define a generating series $J_A$ for the iterated integrals over a rectangle as an element of the Hopf algebra (or the ring). 

Given an open set $U$ in $\R^2$ consider rectangles that lie inside $U$, whose sides are paral
Then $$J_A\times_1 J_B=i(J_{A\times_1 B}).$$

In section $2$ we consider differential $2$ forms on a manifold $X$ and we iterate these forms on a membrane $M$ of dimension $2$. We define a membrane to be the image of a rectangle under a map $\theta$ which is continuous everywhere and smooth except on a finite union of lines such that $\theta$ restricted to any of these lines is smooth except on finitely many points. We show that under some assumptions on the differential forms the iterated integral over the membrane $M$ depends only on the homotopy class of $M$, where the homotopy keeps the boundary points of $M$ fixed.

Define fundamental $2$-groupoid on $X$ $\Pi_2(X)$ to be the set of membranes $M$ in $X$ considered up to homotopy fixing the boundary. Similarly to the compositions $\times_1$ and $\times_2$ that we defined on rectangles, on the fundamental $2$ groupoid we have two compositions. We can compose two membranes, when they have a common face and the parametrization on he face is the same. We have again horizontal and vertical compositions which we denote again by $\times_1$ and $\times_2$. In case that non of the iterated integrals vanishes for a membrane $M$ in $X$ we obtain a fibration of Hopf algebras over the fundamental $2$-groupoid. If we have a vanishing of any of the iterated integrals, then we only obtain a ring fibration over the fundamental $2$-groupoid $\Pi_2(X)$. The products $\times_1$ and $\times_2$ extend to the rings. If $M_1$ and $M_2$ are two membranes with common face in horizontal direction, then again
$$J_{M_1}\times_1 J_{M_2}=i(J_{M_1\times_1 M_2}).$$

In section $3$ we use the products $\times_1$ and $\times_2$ on the rings to show that for $SL_3(\Z)$ the generating series of iterated integrals over a geodesic membrane satisfies a cocycle condition.

In section $4$ we consider admissible membranes on $X(\C)$, where $X$ is an algebraic variety and an admissible membrane means that the boundary of the membrane $M$ lie on the complex points of a divisor of $X$. We show that the iterated integrals can be interpreted as periods of mixed Hodge structures of a relative homology of algebraic varieties over $\bar{\Q}$. 

In section $5$ we give a motivic interpretation of iterated integrals over a membrane of Hilbert modular forms associated to real quadratic fields. We define a non-commutative modular symbol as the generating series of these iterated integrals. 

In section $6$ we give a definition of  multiple completed Dedekind zeta functions over a field $K$ for $K$ equal to $Q$, for $K$ a real quadratic field of class number $1$ and imaginary quadratic field of class number $1$ and $K$ a CM field of class number $1$. We give an interpretation of these multiple completed Dedekind zeta functions as periods of mixed Hodge structures that come from algebraic varieties.

{\bf Acknowledgments.}  
I would like to thank Professor Yuri Manin for the inspiring talk that he gave on non-commutative modular symbol, and for the several conversations we had. I would like to thank also Professor Ronnald Brown for the examples for higher categories without which I would not be able to start this paper. Professor Alexander Goncharov pointed out to me that I should look for a motivic interpretation. I am grateful to him for this suggestion.

This work was initiated at Max-Planck Institute f\"ur Mathematik.
I am very grateful for the stimulating atmosphere, created there,
as well as for the financial support during my stay. Many thanks
are due to the University of Durham for the kind hospitality
during the academic year 2005-2006, when part of this work was
done, and to the Arithmetic Algebraic Geometry Marie Curie Network
for the financial support.


\section{Iteration over a membrane}
\subsection{Definitions and basic properties}
Consider two 2-forms $\a=f(x,y)dx\wedge dy$ and 
$\b=g(x,y)dx\wedge dy$ defined on $(x,y)\in[a_x,b_x]\times[a_y,b_y]$. Let $$\Delta_x=\{(x^1,x^2)\in \R^2|a_x\leq x^1\leq x^2\leq b_x\}$$ and 
$$\Delta_y=\{(y^1,y^2)\in \R^2|a_y\leq y^1\leq y^2\leq b_y\}$$
be two simpleces. 
We define iteration of $\a$ and $\b$ to be
$$\int\dots\int_{\Delta_x\times \Delta_y}
f(x^1,y^1)g(x^2,y^2)dx^1\wedge dy^1\wedge dx^2\wedge dy^2.$$
Similarly to the iterated integrals over a path, we want to express
a product of iterated integrals over a rectangle as a sum of iterated integrals over a rectangle. Consider the simplest case of a product of 
$$\int\int_{
\begin{tabular}{c}$0<x_1<1$\\$0<y_1<1$\end{tabular}}f(x_1,y_1)dx_1\wedge dy_1$$
with
$$\int\int_{
\begin{tabular}{c}$0<x_2<1$\\$0<y_2<1$\end{tabular}}g(x_2,y_2)dx_2\wedge dy_2.$$
This product can be expressed as a sum of four iterated integrals
$$\int\dots\int_{
\begin{tabular}{c}$0<x_1<x_2<1$\\$0<y_1<y_2<1$\end{tabular}}f(x_1,y_1)g(x_2,y_2)dx_1\wedge dy_1\wedge dx_2\wedge dy_2,$$
$$\int\dots\int_{
\begin{tabular}{c}$0<x_1<x_2<1$\\$0<y_2<y_1<1$\end{tabular}}f(x_1,y_1)g(x_2,y_2)dx_1\wedge dy_1\wedge dx_2\wedge dy_2,$$
$$\int\dots\int_{
\begin{tabular}{c}$0<x_2<x_1<1$\\$0<y_1<y_2<1$\end{tabular}}f(x_1,y_1)g(x_2,y_2)dx_1\wedge dy_1\wedge dx_2\wedge dy_2$$
and
$$\int\dots\int_{
\begin{tabular}{c}$0<x_2<x_1<1$\\$0<y_2<y_1<1$\end{tabular}}f(x_1,y_1)g(x_2,y_2)dx_1\wedge dy_1\wedge dx_2\wedge dy_2.$$
Note that the difference between these four integrals is in the domain of integration. We are going to call each of them an iterated integral over a rectangle.
\bde{1.1}
Given $n$ $d$-forms $\a_i(x^i_1,\dots,x^i_d)$ for $i=1,\dots,n$ on $A=[a_1,b_1]\times\dots \times[a_d,b_d]$ and $d$
permutations $\si_1,\dots,\si_d$ of the set $\{1,2,\dots,n\}$ of $n$ elements, let 
$$\Delta(\si_j)=\{(x^1_j,x^2_j,\dots,x^n_j|
a_j\leq x^{\si_j(1)}_j\leq x^{\si_j(2)}_j\leq\dots\leq x^{\si_j(n)}_j\},$$
for $j=1,\dots,d$.
We define an iterated integral over  $[a_1,b_1]\times\dots \times[a_d,b_d]$ to be
$$I_A(\a_1,\dots,\a_n,\si_1,\dots\si_d)=
\int\dots\int_{\Delta(\si_1)\times\dots\times\Delta(\si_d)}
\\a_1(x^1_1,\dots,x^1_d)\dots\wedge\a_i(x^n_1,\dots,x^n_d).$$
\ede
\bde{1.2}
Given two permutations $\si$ of the set $\{1,2,\dots,n_1\}$ of $n_1$ elements and $\tau$ of the set $\{n_1+1,n_1+2,\dots,n_1+n_2\}$ of $n_2$ elements, denote by 
$\Si(\si)(\tau)$ the set of all shuffles which consist
of all permutations $\rho$ such that
if $x_i$'s for $i=1,\dots,n_1+n_2$ satisfy 
$$0<x_{\rho(1)}<\dots <x_{\rho(n_1+n_2)}<1$$
then they satisfy 
$$0<x_{\si(1)}<\dots <x_{\si(n_1)}<1$$
and 
$$0<x_{\tau(n_1+1)}<\dots <x_{\tau(n_1+n_2)}<1.$$
The next proposition is easy to prove but very useful.
\ede
\bpr{1.2} (Shuffle relation)
Let $\a_1\dots\a_{n_1+n_2}$ 
be d-forms on $A=[a_1,b_1]\times\dots\times[a_d,b_d]$. 
Let $\si_1,\dots,\si_d$ be permutations of
$\{1,\dots,n_1\}$, 
and let $\tau_1,\dots,\tau_d$ be permutations of
$\{n_1+1,\dots,n_1+n_2\}$. Let 
$$I_1=I_A(\a_1,\dots,\a_{n_1},\si_1,\dots,\si_d)$$ and
$$I_2=I_A(\a_{n_1+1},\dots,\a_{n_1+n_2},\tau_1,\dots,\tau_d)$$ be two iterated integrals. Then
$$I_1I_2=\sum_{
\rho_j\in\Si(\si_j)(\tau_j),\mbox{ for }j=1,\dots,d
}
I_A(\a_1,\dots,\a_{n_1+n_2},\rho_1,\dots,\rho_d).$$
\epr
\subsection{Composition of rectangles and composition of iterated integrals}
Consider the following composition of rectangles. Let $a_1<b_1<c_1$ and $a_2<b_2<c_2$ be real numbers. We can compose horizontally the rectangle $A=[a_1,b_1]\times[a_2,b_2]$ with 
the rectangle $B=[b_1,c_1]\times[a_2,b_2]$ to obtain $A\times_1B=[a_1,c_1]\times[a_2,b_2]$. Also we can compose vertically the rectangle $A=[a_1,b_1]\times[a_2,b_2]$ with the rectangle $C=[a_1,b_1]\times[b_2,c_2]$ to obtain $A\times_2C=[a_1,b_1]\times[a_2,c_2]$. 
Let $D=[b_1,c_1]\times[b_2,c_2]$. Then
 $(A\times_1 B)\times_2(C\times_1 D)=(A\times_2 C)\times_1(B\times_2 D)$.
This composition of rectangles is an example of compositions of $2$-morphisms in a bi-category. 
We are going to show that iterated integrals over rectangles compose as $2$-morphisms in a 
bi-category. Similarly, if we consider $d$-dimensional parallelepipeds which are 
products of $d$ intervals, we can define $d$ different compositions - one in each direction of a coordinate axis. The composition of two such parallelepipeds can be defined only when they have a common face which is their intersection. This composition can be interpreted as a composition of $d$-morphisms in a cubical $d$-category. We will concentrate our attention to the two-dimensional case which will be done in more details. And we will only mention how this generalizes to higher dimensions.

Consider the integral 
$$I=I_{A\times_1B}(\b_1,\dots,\b_{n},\rho_1,\rho_2).$$
Let 
$$\Delta^i(\rho_1)=
\{(x_1,\dots,x_{n})|a_1\leq x_{\rho_1(1)}\leq \dots \leq x_{\rho_1(i)}\leq b_1\leq
x_{\rho_1(i+1)}\leq\dots\leq x_{\rho_1(n)}\leq c_1\},$$
for $i=1,\dots,n-1$. For $i=0$ and $i=n$ define
$$\Delta^0(\rho_1)=
\{(x_1,\dots,x_{n})|b_1\leq x_{\rho_1(1)}\leq\dots\leq x_{\rho_1(n)}\leq c_1\},$$
and
$$\Delta^{n}(\rho_1)=
\{(x_1,\dots,x_{n})|a_1\leq x_{\rho_1(1)}\leq\dots\leq x_{\rho_1(n)}\leq b_1\}.$$
Let, also,
$$\Delta(\rho_2)=
\{(y_1,\dots,y_{n})|a_2\leq y_{\rho_2(1)}\leq\dots\leq y_{\rho_2(n)}\leq b_2\}.$$
\bde{2.0}
With the above definition of $\Delta^i(\rho_1)$ and $\Delta(\rho_2)$, define
$$I^i_{A\times_1 B}(\b_1,\dots,\b_{n},\rho_1,\rho_2)$$ to be
$$\int\dots\int_{\Delta^i(\rho_1)\times\Delta(\rho_2)}
\b_1(x_1,y_1)\wedge\dots\wedge\b_n(x_n,y_n).$$
\ede
\ble{2.1}
With the above notation we have
$$I_{A\times_1B}(\b_1,\dots,\b_{n},\rho_1,\rho_2)=
\sum_{i=0}^n I^i_{A\times_1 B}(\b_1,\dots,\b_{n},\rho_1,\rho_2).$$
\ele

\bde{2.1}
If $\si$ is a permutation of $m$ elements $\{1,\dots, m\}$ and $\tau$ is a permutation of $n$ elements
$\{1,\dots,n\}$, define $\rho$ to be the unique permutation of $m+n$ elements with the following property: $\rho(i)=\si(i)$ for $i=1,\dots,m$ and 
$\rho(m+i)=\tau(i)$ for $i=1,\dots,n$. We shall denote such permutation $\rho$ by
$(\si,\tau)$. 
\ede

\ble{2.2}
With the above sets $A$ and $B$ we define
$$I_1=I_{A}(\b_1,\dots,\b_{m},\si_1,\si_2)$$ and
$$I_2=I_{B}(\b_{m+1},\dots,\b_{m+n},\tau_1,\tau_2).$$
Then 
$$I_1I_2=\sum_{\rho_2\in\Sigma(\si_2)(\tau_2)} 
I^m_{A\times_1 B}(\b_1,\dots,\b_{m+n},(\si_1,\tau_1),\rho_2).$$
\ele

Now we shall define a generating series for the above type of iterated integrals. Let
$\a_1,\dots,\a_k$ be $k$ distinct $2$-forms defined on $A=[a_1,b_1]\times[a_2,b_2]$.
Consider the non-commutative free polynomial ring 
$$S=\C<z_1,\dots,z_k>,$$ 
where 
$z_i$ corresponds to $\a_i$. Out of this non-commuting ring we shall construct a ring
$R_A$ which captures the algebraic properties of iterated integrals of the $2$-forms
$\a_1,\dots,\a_n$ over the domain $A$. We are going to show also that $R_A$ has the structure of a Hopf algebra.

Each monomial that we are going to consider will have a coefficient $1$ in front of the variables.
For each monomial $M\in S$ of degree $n$, consider $3$ permutations $\si_0$, $\si_1$ and $\si_2$ of $\{1,2,\dots,n\}$. Define an action of $\si_0$ on $M$ which permutes the order of multiplication in the monomial $M$ of degree $n$ in $z_1,\dots,z_k$. Denote this action by $M^{\si_0}$. Define action of $\si_0$ on the right on $\si_1$ and $\si_2$. So that $\si_i$ is mapped to $\si_i \si_0$ for $i=1,2$. Consider the triple
$(M,\si_1,\si_2)$. As we mentioned earlier $z_i$ corresponds to $\a_i$. We shall denote it by $\mu(z_i)=\a_i$ We extend this correspondence to monomials in $S$. For example, if $M=z_2z_1^2z_2$ then we define $\mu(M)=(a_2,a_1,a_1,a_2)$. Denote by
$I_A(M,\si_1,\si_2)$ the iterated integral $I_A(\b_1,\dots, \b_n,\si_1,\si_2)$ 
where $(b_1,\dots,\b_n)=\mu(M)$. Note that $$I_A(M^{\si_0},\si_1\si_0,\si_2\si_0)=I_A(M,\si_1,\si_2).$$ For that reason we define
an equivalence $(M^{\si_0},\si_1\si_0,\si_2\si_0)\sim(M,\si_1,\si_2)$. We shall denote the equivalence class of $(M,\si_1,\si_2)$ by $[M,\si_1,\si_2]$. 
Consider the $\C$-module $V$ generated by $[M,\si_1,\si_2]$. We want to remark that if $M=1$ we do have an element $[M,\si_1,\si_2]$,where  $\si_1$ and $\si_2$ are permutation of $0$ number of elements. 
We shall call 
$[M,\si_1,\si_2]$ a monomial of degree $n$ if the monomial $M$ is of degree $n$. Let $R^0_A$ be a $\Z$-module generated by all 
$$r_A(M,\si_1,\si_2)=I_A(M,\si_1,\si_2)[M,\si_1,\si_2]$$ from $V$,
where $M$ is any monomial in $S$ and $\si_1$ and $\si_2$ are any permutations of $n$ elements with
$n=degree(M)$. If $M$ is of degree $0$ we consider $r(M,\si_1,\si_2)$ to be $1$. 
For technical reasons we assume that non of the iterated integrals vanishes.

We are going to define a ring structure on the $\Z$-module
$R^0_A$. Let $M$ be a monomial in $S$ of degree $m$ and let $\si_1$ and $\si_2$ be two permutations of $\{1,2,\dots,m\}$. Let $N$ be a monomial in $S$ of degree $n$ and let $\tau_1$ and $\tau_2$ be two permutations of $\{1,2,\dots,n\}$. 
We define multiplication of monomial in $R^0_A$ in the following way:
\bde{1.7}$$r(M,\si_1,\si_2)r(N,\tau_1,\tau_2)=
\sum_{\rho_i\in \Si(\si_i)(\tau_i),i=1,2}
r(MN,\rho_1,\rho_2).$$
\ede
This corresponds to multiplication of iterated integrals over $A$ as in Proposition 1.3
This multiplication is associative. One easily checks that
$$r(M,\si_1,\si_2)r(N,\tau_1,\tau_2)r(P,\zeta_1,\zeta_2)=
\sum_{\rho_i\in \Si(\si_i)(\tau_i)(\zeta_i),i=1,2}
r(MNP,\rho_1,\rho_2),$$
where $\Si(\si_i)(\tau_i)(\zeta_i)$ is a triple shuffle of the permutations 
$\si_i$, $\tau_i$ and $\zeta_i$.

Let $I(n)$ be the ideal in $R^0_A$ generated by all monomials $r\in R^0_A$ of degree greater than $n$. We have $I(n+1)\subset I(n)$ for all $n$. Define $R_A$ to be the inverse limit of $R^0_A/I(n)$ when $n$ tends to infinity. In $R_A$ we shall denote the monomials again by $r(M,\si_1,\si_2)$. If $M=1$ then we set $r(M,\si_1,\si_2)=1$.

Now we are ready to define generating series for iterated integrals over a rectangle.
Let 
$$J_A(\a_1,\dots,\a_k)=
\sum_{r_A(M,\si_1,\si_2)\in R_A}
r(M,\si_1,\si_2).$$
The sum converges in the induced topology on $R_A$ from the inverse limit of the discrete rings $R^0_A/I(n)$.
\bth{1.8}
The ring $R_A$ has the structure of a Hopf algebra.
\eth
\proof The monomial $r_A(M,\rho_1,\rho_2)$ is the same as $r_A(M^{\rho_1^{-1}},(1),\rho_2\rho_1^{-1})$. We are going to define coproduct on
$r_A(M,(1),\rho_2)$. Let $r_A(M,(1),\rho_2)$ be of degree $n$. For $i=1,\dots,n-1$ define
$$\Delta^i(r_A(M,(1),\rho_2))=
r_A(M',(1),\si_2)\otimes r_A(M'',(1),\tau_2),$$
where $degree(M')=i$, $\si_2$ is a permutation of $i$ elements, $degree(M'')=n-i$,
$\tau_2$ is a permutation of $N-i$ elements, and $\rho_2=(\si_2,\tau_2)$. Note that if such $\si_2$ and $\tau_2$ exist then they are unique. If they do not exist we define
$$\Delta^i(r_A(M,(1),\rho_2))=0.$$
Define also $$\Delta^0(r_A(M,(1),\rho_2))=1 \otimes r_A(M,(1),\rho_2),$$
and
$$\Delta^n(r_A(M,(1),\rho_2))=r_A(M,(1),\rho_2)\otimes 1,$$
Define the coproduct in $R_A$ to be
$$\Delta=\sum_{i=0}^n \Delta^i$$ for a monomial of degree $n$.
It is easy to check that the coproduct is coassociative. 
Define the counit $$\varepsilon:R_A \rightarrow \Z$$ to be zero on monomials of degree greater than $0$ and $\varepsilon(1)=1$.
It is straight forward to check that $\varepsilon\otimes id\circ\Delta=id$.
So $R_A$ is a bialgebra.

Now we define the antipode $S:R_A \rightarrow R_A$. We call an monomial $r$ Lie like if
$\Delta(r)=1\otimes r +r \otimes 1$. For Lie like elements we define $S(r)=-r$. 
It is straight forward to check that for Lie like elements $r$
$$m\circ S\otimes id \circ \Delta(r)=0=u\circ\varepsilon(r),$$
where $m$ is the multiplication in the ring and $u:Z\rightarrow R_A$ the unique homomorphism. To define $S$ for all monomials we use induction on the number of summands of $\Delta(r)$ for arbitrary monomial $r\in R_A$. From 
$$m\circ S\otimes id \circ \Delta(r)=0$$ we can express $S(r)$ is terms of $S(r')$
where $\Delta(r')$ has less summands that $\Delta(r)$. It is clear that $S$ satisfies the axioms for a antipode. Thus, $R_A$ is a commutative Hopf algebra.
\bco{1.10}
The genrating series $J_A(\a_1,\dots,\a_k)$is a group-like elements in the Hopf algebra. That is,
$$\Delta J_A(\a_1,\dots,\a_k=J_A(\a_1,\dots,\a_k)\otimes J_A(\a_1,\dots,\a_k).$$
\eco
The proof is straight forward.

We can compose horizontally the rectangle $A=[a_1,b_1]\times[a_2,b_2]$ with 
the rectangle $B=[b_1,c_1]\times[a_2,b_2]$ to obtain $A\times_1B=[a_1,c_1]\times[a_2,b_2]$. We are going to relate $R_A$, $R_B$ and $R_{A\times_1 B}$. For that reason we are going to define a new Hopf algebra
$R_{A,B}$ in terms of iterated integrals of the type 
$$I^i_{A\times_1 B}(\b_1,\dots,\b_{n},\si_1,\si_2).$$
Consider the $\C$-module $V$ generated by $[M,\si_1,\si_2]$ and an index $i$, where $i$ takes values in the set $\{0,1\dots,n\}$ and $n=degree(M)$. Let
$$r^i_{A,B}(M,\si_1,\si_2)=I^i_{A\times_1 B}(\mu(M),\si_1,\si_2)$$
Let $R^0_{A,B}$ be the $\Z$-submodule of $V$ generated by $r^i_{A,B}(M,\si_1,\si_2)$.
We shall call $r^i_{A,B}(M,\si_1,\si_2)$ a monomial of degree $n$ is $M$ is of degree $n$. 

We would like to define multiplication in $R^0_{A,B}$. In order to do that we need some definitions. Let $\si_1$ and $\tau_1$  be permutations of  of the sets $\{1,\dots i,i+j+1\dots,m+j\}$ and $\{i+1,\dots,i+j,\dots,m+j+1,m+n\}$, respectively. Define restricted suffle, which we denote by $\Sigma(\si_1)(\tau_1)\{i,j\}$, the set of permutation $\rho_1$ such that 
$$x_{\rho_1(1)}<\dots<x_{\rho_1(i+j)}<a<x_{\rho_1(i+j+1)}<\dots<x_{\rho_1(m+n)}$$ implies
$$x_{\si_1(1)}<\dots<x_{\si_1(i)}<a<x_{\si_1(i+j+1)}<\dots<x_{\si_1(m+j)}$$ and 
$$x_{\rho_1(i+1)}<\dots<x_{\rho_1(i+j)}<a<x_{\tau_1(m+j+1)}<\dots<x_{\tau_1(m+n)}.$$
\bde{1.10}We define multiplication in $R^0_{A,B}$ 
$$r^i_{A,B}(M,\si_1,\si_2)r^j_{A,B}(N,\tau_1,\tau_2)=
\sum_{\rho_1\in \Si(\si_1)(\tau_1)\{i,j\},\rho_2\in \Si(\si_2)(\tau_2)}
r^{i+j}_{A,B}((MN)^{\rho_0},\rho_1,\rho_2),$$
Where $rho_0$ permutes the variables of $MN$ so that the first $i$ number of variables 
after the permutation are the first $i$ variables of $M$, the next $j$ number of variables after the permutation are the first $j$ number of variables of $N$, the next $m-i$ number of variables after the permutation are the last $m-i$ number of variables  in $M$ , and finally, the last $n-j$ number of variables after the permutation are the last $n-j$ elements of $N$.
\ede 
Let $I(n)$ be the ideal in $R^0_{A,B}$ generated by monomials of degree greater that $n$. Let $R_{A,B}$ be the ring obtained by taking the inverse limit of $R^0_{A,B}/I(n)$. 
\bpr{1.11}
The ring $R_{A,B}$ has the structure of a Hopf algebra.
\epr
We omit the proof since it is essentially the same structure as the Hopf algebra $R_{A\times_1 B}$.

We want to define a multiplication so that when we multiply $R_A$ and $R_B$ we obtain
$R_{A\times_1 B}$. This is not possible but something very close to it can be done.
We can define homomorphism of Hopf algebras  
$$\times_1:R_A \otimes R_B \rightarrow R_{A,B}$$ 
and
$$i:R_{A\times_1 B}\rightarrow R_{A,B},$$
so that 
$$J_A(\a_1,\dots,\a_k)\times_1 J_B(\a_1,\dots,\a_k)=
i(J_{A\times_1 B}(\a_1,\dots,\a_k)).$$

We shall consider the following coproduct on $R_A\otimes R_B$. We are going to use the notation from the proof of Theorem 1.9. Let $r_A$ and  $r_B$ be monomials in $R_A$ and $R_B$, respectively. Let $n=degree(r_B)$.  Let $\Delta^i_A$ be the $i$-th component of the coproduct on $R_A$ and similarly, $\Delta^j_B$ be the $j$-th component of the coproduct on $R_B$. The coproduct on $R_A\otimes R_B$ that we are going to consider is 
$$\Delta=\sum_{i=0\mbox{ or }j=n}\Delta^i_A\otimes \Delta^j_B.$$
The Lie-like elements in $R_A\otimes R_B$ are the ones coming from the embeddings of $R_A$ and $R_B$.
It is easy to check that the maps $R_A \rightarrow R_A\otimes R_B$ and 
$R_B \rightarrow R_A\otimes R_B$ are homomorphisms of Hopf algebras.

Consider the homomorphism $$\times_1:R_A \otimes R_B \rightarrow R_{A,B}$$ 
defined by 
$$r_A(M,(1),\si_2)\otimes r_B(N,(1),\tau_2)=
\sum_{\rho_2\in\Sigma(\si_2)(\tau_2)}r^m_{A,B}(MN,(1),\rho_2),$$
where $m=degree(M)$. This definition follows closely the property of iterated integrals described in Lemma 1.7.

Consider the map
$$i:R_{A\times_1 B}\rightarrow R_{A,B}$$
defined on monomials by 
$$i(r_{A\times_1 B}(M,(1),\si_2)=\sum_j r^j_{A,B}(M,(1),\si_2).$$
It is easy to check that this is a homomorphism of Hopf algebras.
\bth{1.5}
Let $A=[a_1,b_1]\times[a_2,b_2]$, $B=[b_1,c_1]\times[a_2,b_2]$ and $A\times_1B=[a_1,c_1]\times[a_2,b_2]$, where
where $a_1<b_1<c_1$ and $a_2<b_2$.
And let $\a_1,\dots,\a_k$ be $2$-forms on $A\times_1B$.
Then 
$$i(J_{A\times_1 B}(\a_1,\dots,\a_k))=
J_A(\a_1,\dots,\a_k)\times_1 J_B(\a_1,\dots,\a_k)$$
\eth
\proof
We are going to show that 
$$r^i_{A, B}(P,\rho_1,\rho_2)$$
occurs exactly once among the monomials in the product
$$J_A(\a_1,\dots,\a_k)J_B(\a_1,\dots,\a_k).$$ Also,each monomial in the product gives rise to exactly one element of the type $r^i_{A,B}$. 

Let $r_A(M,\si_1,\si_2)$ and $r_B(N,\tau_1,\tau_2)$ be two monomials in $R_A$ and $R_B$, of degree $m$ an $n$, respectively. We have
$$r_A(M,\si_1,\si_2)r_B(N,\tau_1,\tau_2)=
\sum_{\rho_1=(\si_1,\tau_1),\rho_2\in\Si(\si_2)(\tau_2)}
r^m_{A, B}(MN,\rho_1,\rho_2).$$ 
We can vary the permutations to show that $r^i_{A, B}(P,\rho_1,\rho_2)$ occurs exactly one in the product. Form the shuffles we have 
$$\left(\frac{(m+n)!}{m!n!}\right)^2$$
variation. From the permutations $\si_1$ and $\si_2$ we have $(m!)^2$. And from
$\tau_1$ and $\tau_2$ we have $(n!)^2$ permutations. When we multiply all of them we obtain all permutations, which are $((m+n)!)^2$. So we have counted each of the permutations.

By abuse of notation we are going to write 
$$J_{A\times_1 B}(\a_1,\dots,\a_k)=
J_A(\a_1,\dots,\a_k)\times_1 J_B(\a_1,\dots,\a_k).$$
Similarly, we define $\times_2$ for vertical compositions of rectangles and for the generating series over rectangles with a common face in vertical direction. For rectangles we have
$$(A\times_1 B)\times_2(C\times_1 D)=(A\times_2 C)\times_1(B \times_2 D).$$
Similarly, for the generating series of iterated integrals over rectangles we have
$$(J_A\times_1 J_B)\times_2(J_C\times_1 J_D)=
(J_A\times_2 J_C)\times_1(J_B \times_2 J_D),$$
where $J_A=J_A(\a_1,\dots,\a_k)$.
\section{Homotopy invariance of iterated integrals over a membrane}
Let $X$ be a manifold. And let 
$\omega_1,\dots\omega_n$ be closed $2$-forms on $X$. Denote
by $A$ the rectangle $[a_1,a_2]\times[b_1,b_2]$ in $\R^2$ for 
$a_1<a_2$ and $b_1<b_2$. Let $\varphi:A\rightarrow X$ be continuous, almost everywhere differentiable function. Let $\a_i=\varphi^*\omega_i$ be the pull back of the form $\omega_i$ to $A$. We can iterated the forms $\omega_1,\dots\omega_n$ on $f(A)$ by pulling them back to $A$.
Let $I(\varphi,A,\omega_1,\dots\omega_n,\si_x,\si_y)$ be the iterated integral defined by $I(A,\a_1,\dots\a_n,\si_x,\si_y)$
Let $\omega_i=f_i \omega$ where $f_i\in A^0(X)$, and $\omega\in A^2(X)$ Let $\varphi_t:A\rightarrow X$ for $t\in[0,1]$ be homotopy of maps from $A$ to $X$ such that the restriction to the boundary does not vary with $t$. 
$$\varphi_t|_{\partial A}=\varphi_0|_{\partial A}.$$ That is, the homotopy is
constant on each boundary point of $A$. 
\bpr{3.1} The iterated integral is homotopy invariant with
respect to the membrane of integration when the homotopy is constant on the boundary. With the above notation
$$I(\varphi_0,A,\omega_1,\dots\omega_n,\si_x,\si_y)=
I(\varphi_1,A,\omega_1,\dots\omega_n,\si_x,\si_y)$$
\epr
\proof Let $A^n$ be a product of $n$ copies of $A$. Let
$\varphi_t^n:A^n\rightarrow X^n$ be $n$ copies of the map $\varphi_t$. Let $\omega_i^i$ be the form $\omega_i$ on $X$ pulled back to the $i$-th copy of $X^n$. 
Let $$\Delta=\{(x_1,y_1,\dots,x_n,y_n): a_1\leq x_{\si_x(1)}\leq \dots x_{\si_x(n)}\leq a_2\mbox{ and }b_1\leq y_{\si_y(1)}\leq \dots y_{\si_y(n)}\leq b_2\}.$$

Then the iterated integral can be written as $$I(\varphi_t,A,\omega_1,\dots\omega_n,\si_x,\si_y)=
\int_{\Delta}(\varphi_t^n)^*\Omega,$$ 
where 
$$\Omega=\omega_1^1\wedge\dots\wedge\omega_n^n$$
This integral is also equal to 
$$\int_{\varphi_t^n(\Delta)}\Omega.$$
Let $\Delta'=\cup_{t\in[0,1]}\varphi_t^n(\Delta)$
In order to show that 
$$\int_{\varphi_t^0(\Delta)}\Omega=
\int_{\varphi_t^1(\Delta)}\Omega,$$
we need to show that $$\int_{\Delta'}\Omega=0.$$
The set $\Delta'$ is obtained by setting one of the inequalities
in the definition of $\Delta$ to an equality, and then map it with $\varphi_t^n$. 

Consider $x_{\si_x(i)}=x_{\si_x(i+1)}$ in $\Delta$ and map it with all $\varphi_t^n$ to $\Delta'$ Call the image $\Delta'_{x,i}.$ Let $$\partial/\partial x_{\si_x(i)}$$
be a directional derivative on the $\si_x(i)$-th factor of $A$
in $A^n$. Let $$\Phi:A^n\times[0,1]\rightarrow X^n$$ be the 
homotopy of $\varphi_t^n$. Let
$$\frac{\partial}{\partial z}=
\Phi_*(\frac{\partial}{\partial x_{\si_x(i)}}),$$
on $\Delta'_{x,i}$. Restrict $\Omega$ to $\Delta'_{x,i}.$ Then locally 
$\omega_{\si_x(i)}|_{\Delta'_{x,i}}=\b_i\wedge dz$ and
$\omega_{\si_x(i+1)}|_{\Delta'_{x,i}}=\b_{i+1}\wedge dz,$ where $\b_i$ and $\b_{i+1}$ are $1$-forms. We have repetition of $dz$ because $x_{\si_x(i)}=x_{\si_x(i+1)}$. Thus $\Omega$ restricted
 to $\Delta'_{x,i}$ is zero because we have wedge of two $dz$'s.
 
We have one more type of boundary component. It is obtained when we set $x_{\si_x(1)}=a_1$, or $x_{\si_x(n)}=a_2$, or
$y_{\si_y(1)}=b_1$, or $y_{\si_y(n)}=b_2$. In any of these cases, we restrict $\Omega$ to the corresponding boundary component. Let $\Delta'_{x,0}$ be the boundary component corresponding  to $x_{\si_x(1)}=a_1$. Since one of the variables is set to a constant, we have that $\Omega|_{\Delta'_{x,0}}=0$. Therefore all boundary components of $\tilde{\Delta}$ give zero except $\varphi_0^n(A^n)$ and
$\varphi_1^n(A^n)$ which give the two iterated integrals. 
By assumption each form $\omega_i$ is closed. Then $d\Omega=0$. So by Stoke's theorem we have
$$0=\int_{\tilde{\Delta}}d\Omega=
\sum_{i=0}^1
(-1)^i I(\varphi_i,A,\omega_1,\dots\omega_n,\si_x,\si_y).$$

\bre{2.5}
Given a membrane $M$ on $X$ which is the image of a rectangle mapped to $X$ we can consider the generating series of iterated integrals $J_M(\omega_1,\dots,\omega_k)$.
If $N$ is another membrane such that $M$ and $N$ have a common face parametrized in the same way then we can 'add' $M$ and $N$ to obtain $M\times_1 N$. For the generating series we obtain 
$$J_M (\omega_1,\dots,\omega_k)\times_1 J_N(\omega_1,\dots,\omega_k) =
J_{M\times_1 N} (\omega_1,\dots,\omega_k).$$
All these generating series depend on the homotopy class of the membranes. Let the fundamental $2$-groupoid on $X$ consist of homotopy maps of rectangles to $X$. We add two elements of the fundamental $2$-groupoid when they have a common face, parametrized in the same way. Thus, 
$$J(\omega_1,\dots,\omega_k)$$ (the subscript is not missing)
is a Hopf algebra fibration over the fundamental $2$-groupoid (or a ring fibration in case of vanishing of an iterated integral). We are going to use that, when we consider the arithmetic groups $GL_m(\Z)$ and in more details $GL_3(\Z)$, in the next section.
\ere
\section{Non-commutative periods for $GL_m$ over the integers and over imaginary quadratic rings}
In this section we will construct non-commutative modular symbols for arithmetic groups $\Gamma$ which are commensurable to $GL_m(\Z)$ or commensurable to $GL_m({\cal{O}}_K)$ where ${\cal{O}}_K$ is the ring of integers in an imaginary quadratic field $K$.

Let ${\cal{H}}_{\R}^m$ be the space of positive definite symmetric $m\times m$ matrices with
real coefficients modulo multiplication by a positive scalar. One can identify 
${\cal{H}}_{\R}^m$ with $SL_m(\R)/SO_m(\R)$. $GL(\R)$ acts on ${\cal{H}}_{\R}^m$ in the following way. Let $H\in{\cal{H}}_{\R}^m$ and $g\in GL(\R)$. Then $$g(H)=g H g^T,$$
where $g^T$ is the transposition of the matrix $g$. 

We shall add points at infinity to ${\cal{H}}_{\R}^m$. 
Let $g\in Mat_{m\times m}(\Q)$ be a matrix with rational coefficients of rank $r=1$,
For such $g$ we construct the semi-definite symmetric $m\times m$ matrix $gg^T$. If $g$ is of rank $1$ then $g=v_1v_2^T$, where $v_i$ $i=1,2$ are $m$ dimensional non-zero vectors. Consider $gg^T=v_1v_2^Tv_2v_1^T=(v_2^Tv_2)v_1v_1^T$. Modulo a scalar multiplication by a positive real, we have that $gg^T=v_1v_1^T.$ Given $n$ non-zero vectors $v_i$ for $i=1,\dots,n$, let $H_i=v_iv_i^T$.
Let $\tilde{{\cal{H}}}_{\R}^m$ be the space of semi-definite symmetric $m\times m$ matrices. Then for any $t_i\in[0,1]$ such that $\sum_i t_i=1$ we have that 
$\sum_i t_iH_i\in \tilde{{\cal{H}}}_{\R}^m$. Let $\overline{{\cal{H}}}_{\Q}^m$
consists of all point $H=\sum_i t_iH_i$ with the above $t_i$'s and $H_i$'s modulo a scalar multiplication by a positive real number. 
Let $$\overline{{\cal{H}}}_{\R}^m=\overline{{\cal{H}}}_{\Q}^m \cup {\cal{H}}_{\R}^m.$$
Let $P^{m-1}(\Q)$ be $(n-1)$-dimensional projective space over the rational numbers.
Identify $P^{m-1}(\Q)$ with a subset of $\overline{{\cal{H}}}_{\R}^m$ in the following way. Let $v$ be  vector in $\Q^n-\{O\}$ considered as a representative of $P^{m-1}(\Q)$. Consider the symmetric matrix $vv^T$ as an element of  
$\overline{{\cal{H}}}_{\Q}^m\subset \overline{{\cal{H}}}_{\R}^m$. 

Let $P_0,\dots,P_k$ be points of $P^{m-1}(\Q)\subset \overline{{\cal{H}}}_{\R}^m$.
Consider the triangle consisting of points $P=\sum_{i=0}^k t_k P_k$ satisfying 
$0\leq t_i \leq 1$ and $\sum_{i=0}^k t_i=0$. Denote this simplex by $\Delta$. 
Denote by $[0,1]^k$ the $k$-dimensional cube. Let $(x_1,\dots,x_k)$ be a point of the $k$-dimensional cube. Consider the map $$\a_k:[0,1]^k \rightarrow \Delta$$ defined by
$$\a_k(x_1,\dots,x_k)=\sum_{i=0}^k t_k P_k,$$
where
$$\begin{tabular}{l}
$t_0=x_1,$\\
$t_1=(1-x_1)x_2,$\\
$t_2=(1-x_1)(1-x_2)x_3,$\\
$\dots$\\
$t_{k-1}=(1-x_1)\dots(1-x_{k-1})x_k,$\\  
$t_k=(1-x_1)\dots(1-x_k).$\\
\end{tabular}$$
The map $\a_k$ is a homeomorphism between the interior of $[0,1]^k$ and $\Delta$. Also each face of the cube is either mapped to a simplex of smaller dimension, or it is homeomorphic to a $k-1$ dimensional face of $\Delta$ with a map $\a_{k-1}$.

Let $K$ bean imaginary quadratic field. And let ${\cal{O}}_K$ be the ring of integers
in $K$.
Denote by ${\cal{H}}_{\C}^m$ be the space of positive definite Hermitian $m\times m$ matrices modulo multiplication by a positive real number. One can identify 
${\cal{H}}_{\C}^m$ with $SL_m(\C)/SU_m(\R)$. $GL(\C)$ acts on ${\cal{H}}_{\C}^m$ in the following way. Let $H\in{\cal{H}}_{\C}^m$ and $g\in GL(\C)$. Then $$g(H)=g H \bar{g}^T,$$
where $\bar{g}^T$ is the complex conjugation of the transposition of the matrix $g$. 

We shall add points at infinity to ${\cal{H}}_{\C}^m$. 
Let $g\in Mat_{m\times m}(K)$ be a matrix with rational coefficients of rank $r=1$,
For such $g$ we construct the semi-definite Hermitian $m\times m$ matrix $g\bar{g}^T$. If $g$ is of rank $1$ then $g=v_1\bar{v}_2^T$, where $v_i$ $i=1,2$ are $m$ dimensional non-zero vectors. Consider $g\bar{g}^T=v_1\bar{v}_2^Tv_2\bar{v}_1^T=(\bar{v}_2^Tv_2)v_1v_1^T$. Modulo a scalar multiplication by a positive real, we have that $g\bar{g}^T=v_1\bar{v}_1^T.$ Given $n$ non-zero vectors $v_i$ for $i=1,\dots,n$, let $H_i=v_iv_i^T$.
Let $\tilde{{\cal{H}}}_{\R}^m$ be the space of semi-definite Hermitian $m\times m$ matrices. Then for any $t_i\in[0,1]$ such that $\sum_i t_i=1$ we have that 
$\sum_i t_iH_i\in \tilde{{\cal{H}}}_{\C}^m$. Let $\overline{{\cal{H}}}_{\K}^m$
consists of all point $H=\sum_i t_iH_i$ with the above $t_i$'s and $H_i$'s modulo a scalar multiplication by a positive real number. 
Let $$\overline{{\cal{H}}}_{\C}^m=\overline{{\cal{H}}}_{K}^m \cup {\cal{H}}_{\C}^m.$$
Let $P^{m-1}(K)$ be $(n-1)$-dimensional projective space over the rational numbers.
Identify $P^{m-1}(K)$ with a subset of $\overline{{\cal{H}}}_{\C}^m$ in the following way. Let $v$ be  vector in $K^n-\{O\}$ considered as a representative of $P^{m-1}(K)$. Consider the symmetric matrix $v\bar{v}^T$ as an element of  
$\overline{{\cal{H}}}_{K}^m\subset \overline{{\cal{H}}}_{\C}^m$. 

Let $P_0,\dots,P_k$ be points of $P^{m-1}(K)\subset \overline{{\cal{H}}}_{\C}^m$.
Consider the triangle consisting of points $P=\sum_{i=0}^k t_k P_k$ satisfying 
$0\leq t_i \leq 1$ and $\sum_{i=0}^k t_i=0$. Denote this simplex by $\Delta(P_0,\dots,P_k)$. 
Denote by $[0,1]^k$ the $k$-dimensional cube. Let $(x_1,\dots,x_k)$ be a point of the $k$-dimensional cube. Consider the map $$\a_k:[0,1]^k \rightarrow \Delta$$ defined by
$$\a_k(x_1,\dots,x_k)=\sum_{i=0}^k t_k P_k,$$
where
$$\begin{tabular}{l}
$t_0=x_1,$\\
$t_1=(1-x_1)x_2,$\\
$t_2=(1-x_1)(1-x_2)x_3,$\\
$\dots$\\
$t_{k-1}=(1-x_1)\dots(1-x_{k-1})x_k,$\\  
$t_k=(1-x_1)\dots(1-x_k).$\\
\end{tabular}$$
The map $\a_k$ is a homeomorphism between the interior of $[0,1]^k$ and $\Delta$. Also each face of the cube is either mapped to a simplex of smaller dimension, or it is homeomorphic to a $k-1$ dimensional face of $\Delta$ with a map $\a_{k-1}$.

Let $\Ga$ be $GL_m(\Z)$ or $GL_m({\cal{O}}_K)$, where ${\cal{O}}_K$ is the ring of integers in an imaginary quadratic field $K$. Let $\Ga_0$ be a finite index subgroup of 
$\Ga$. Let $\overline{{\cal{H}}}^m$  be $\overline{{\cal{H}}}_{\R}^m$ or $\overline{{\cal{H}}}_{\C}^m$ with boundary points $P^{m-1}(K)$ where $K$ is $\Q$ or an imaginary quadratic field. Let $P\in P^{m-1}(K)\subset \overline{{\cal{H}}}^m$.
Let $g_i\in\Ga$ for $i=0,1,2$. Let $\Delta=\Delta(g_0 P,g_1 P,g_2 P)$. Let $\omega$ be a smooth $2$-forms on $\overline{{\cal{H}}}^m$ which is cusp form with respect to 
$\Ga_0$. (it might represent a zeroth cohomological class in $H^2(\Ga_0,V)$.) Let 
$\omega_i=f_i \omega$ for $i=1,\dots,n$, where $f_i$ is a smooth function on $\overline{{\cal{H}}}^m$, so that $\omega_i$ is a cusp form with respect to $\Ga_0$.
Consider the generating series for the iterated integrals over $\Delta$ of the forms
$\omega_i$ for $i=1,\dots,n$. Let
$$J_P(g_0,g_1,g_2)=J_{\Delta}(\omega_1,\dots,\omega_n).$$
\bth{5.1}
The generating series $J_P$ is a second non-abelian cohomological class in 
$H^2(\Ga,R),$ where $R$ consists of Hopf algebras fibered over the fundamental $2$-groupoid of $\overline{{\cal{H}}}^m$ with compatible multiplications (as defined in sections $2$ and $3$).
\eth
\proof We are not going to use the Hopf algebra structure now. However, we are going to use the two multiplication $\times_1$ and $\times_2$ in the fundamental $2$-groupoid. These two multiplications induce two multiplications on the generating series of iterated integrals. We denote the two induced multiplications again
by $\times_1$ and $\times_2$.

If we have two membranes $M_1$ and $M_2$ with a common face in horizontal direction, we can "add" them to obtain $M_1\times_1 M_2$. Then we have
$$J_{M_1}(\omega_1,\dots,\omega_n)\times_1 J_{M_2}(\omega_1,\dots,\omega_n) 
=J_{M_1\times_1 M_2}(\omega_1,\dots,\omega_n).$$
Similarly, if we have two membranes $M_3$ and $M_4$ with a common face in vertical direction, we can "add" them to obtain $M_3\times_2 M_4$. Then we have
$$J_{M_3}(\omega_1,\dots,\omega_n)\times_1 J_{M_4}(\omega_1,\dots,\omega_n) 
=J_{M_3\times_2 M_4}(\omega_1,\dots,\omega_n).$$
To show that $J_P$ satisfies a $2$-cocycle condition, take four elements $g_i$ for $i=0,1,2,3$ of the group $\Ga$. Let 
$$\Delta_0=\Delta(g_1 P,g_2 P,g_3 P),$$
$$\Delta_1=\Delta(g_0 P,g_2 P,g_3 P),$$
$$\Delta_2=\Delta(g_0 P,g_1 P,g_3 P),$$ and
$$\Delta_3=\Delta(g_0 P,g_1 P,g_2 P).$$
We have that 
$$J_P(g_0,\dots,\hat{g}_i,\dots,g_3)=J_{\Delta_i}(\omega_1,\dots,\omega_n),$$
for $i=0,1,2,3$.
We have to prove an identity of the type
$$J_P(g_0,g_1,g_2)\mbox{ "+" }J_P(g_0,g_2,g_3)=
J_P(g_0,g_1,g_3)\mbox{ "+" }J_P(g_1,g_2,g_3).$$
We are going to prove that
$$J_P(g_0,g_1,g_2)\times_1 J_P(g_0,g_2,g_3)=
J_P(g_0,g_1,g_3)\times_2 J_P(g_1,g_2,g_3).$$
This is equivalent to
$$J_{\Delta_3}(\omega_1,\dots,\omega_n)\times_1 J_{\Delta_1}(\omega_1,\dots,\omega_n)=
J_{\Delta_2}(\omega_1,\dots,\omega_n)\times_2 J_{\Delta_0}(\omega_1,\dots,\omega_n)$$
In other words, we claim that
$$J_{\Delta_3\times_1 \Delta_1}(\omega_1,\dots,\omega_n)=
J_{\Delta_2\times_2 \Delta_0}(\omega_1,\dots,\omega_n),$$
which is true because
$$\Delta_3\times_1 \Delta_1$$ is homotopic to $$\Delta_2\times_2 \Delta_0$$ with homotopy which leaves the boundary fixed. This follows simply from the definition of the map $\a_2$. This finishes the proof.
\bth{5.6}
The cohomological class of $J_P$ is independent of the choice of $P$. In other words, if $Q$ is a boundary point then $J_P$ differs from $J_Q$ by a coboundary.
\eth
\proof
Let $g_0$ and $g_1$ be two elements of $\Ga$.
Let $$\Delta(g_0,g_1)=\Delta(g_0 Q, g_0 P, g_1 P)$$ and 
$$\Delta'(g_0,g_1)=\Delta(g_1 P, g_0 Q, g_1 Q).$$
Define 
$$K_{P,Q}(g_0,g_1)=
J_{\Delta(g_0,g_1)}(\omega_1,\dots,\omega_n)\times_1 J_{\Delta'(g_0,g_1)}(\omega_1,\dots,\omega_n).$$
If $g_0,g_1$ and $g_2$ are three elements of $\Ga$, we can define
$$\partial K_{P,Q}(g_0,g_1,g_2)=
J_{\Delta(g_0,g_1,g_2)}(\omega_1,\dots,\omega_n),$$
where $$\Delta(g_0,g_1,g_2)=
\Delta(g_0,g_1)\times_1 \Delta'(g_0,g_1)\times_1
\Delta(g_1,g_2)\times_1 \Delta'(g_1,g_2)\times_1
\Delta(g_2,g_0)\times_1 \Delta'(g_2,g_1).$$
When we multiply $\partial K_{P,Q}(g_0,g_1,g_2)$ with $J_P(g_0,g_1,g_2)$ and 
$J_Q(g_0,g_1,g_2)$, we are going to integrate over a contractable space with no boundary. So the integral will be zero. In order to make this argument, we need some combinatorics for multiplying many simpleces. To proof Theorem 4.2 the combinatorics is not very difficult. For a general arrangement we have the following lemma.
\ble{5.7}
Suppose we are given elements $M_i$ for $i=1,\dots,n$ of the fundamental $2$-groupoid on $X$. Suppose that each face of $M_i$ is either a point or there exist unique $j$ such that $M_i$ and $M_j$ have a common face and the parametrization of that face coming from $M_i$ and $M_j$ coincides. Then there exists an element $M$ of the fundamental $2$-groupoid which geometrically is a union of $M_i$'s such that the parametrization of $M$ can be subdivided into $M_i$'s. If further more the element $M$ is contractible and simply connected then the product of the iterated integrals $J_{M_i}$ will be $J_M=0$. 
\ele

\bre{5.8}
For $\Ga$ a finite index subgroup of $GL(3,\Z)$ and a self-dual representation $V$
Consider the cuspidal cohomology $H^2_{cusp}(\Ga,V)$. Let $\a_1,\dots,\a_k$ be $2$-forms representing cups cohomological classes such that $\a_i=f_i\a$, where $f_i$ is a smooth form on $SL_3(\R)/SO_3(\R)$ and $\a$ is a smooth $2$-form on $SL_3(\R)/SO_3(\R)$. Them we propose $J_P(\a_1,\dots,\a_k)$ to be the non-abelian modular symbol for $GL_{3/\Q}$.
\ere
\section{Motivic interpretation of iteration over a membrane}
In this section we give a construction of motivic version of iteration over a membrane. 

Let $X$ be a quasi-projective variety of dimension $n$ over $\overline{\Q}$. All the constructions of new varieties will be algebraic. So that all the new varieties can be defined over $\overline{\Q}$. Let $M$ be an admissible  membrane on $X$ of real dimension $n$. We are going to define what admissible membrane means. Let $M$ is the image of 
$\theta:[0,1]^n \rightarrow X(\C)$ such that the restriction of $\theta$ to any of the faces of the cube lies on the complex points of a union of subvarieties of $X$ of codimension $1$. Let 
$$s^i_0=\{(t_1,\dots,t_n)|t_i=0, \mbox{ and }t_j\in[0,1]\mbox{ for }j\neq i\}$$
and
$$s^i_1=\{(t_1,\dots,t_n)|t_i=1, \mbox{ and }t_j\in[0,1]\mbox{ for }j\neq i\}.$$
Let $$Z^i_0\supset{\theta(s^i_0)}$$ and
$$Z^i_1\supset{\theta(s^i_1)},$$
where $Z^i_0$ and $Z^i_1$ are union of varieties of codimension $1$ in $X$ for all $i=1,\dots n$. 
These are the conditions that an admissible membrane $M$ satisfies. 

We want to mention
that iteration over $M$ of $n$-forms from $\Ga(\Omega^n_X,X)$ on $X(\C)$ gives the same value for the integral as iteration over $M'$ if $M$ is homotopic to $M'$ so that the boundary stays fixed in the homotopy. The invariance of the iterated integrals goes even further. If we iterate
$n$-forms on $X(\C)$ over $M$, we can take a homotopy  $M_t$ so the faces 
of $M_t$ lies on the corresponding $Z^i_0$ or $Z^i_1$ for faces of $M_t$ for any $t\in[0,1]$ with $M_0=M$. Then for any $t\in[0,1]$ the iterated integral will have the same value,
because if we integrate $n$-forms from $\Ga(\Omega^n_X,X)$ on a membrane that lies inside a variety of codimension $1$ in $X$ we obtain $0$. 

Consider a membrane $M$ of the above type. For technical reasons, also for the purpose of our applications, we assume that no pair of divisors among $Z^i_0$ and $Z^i_1$ for all $i=1,\dots n$ have a common component. Enlarge any of the divisors if necessary so that for each $i$ there is a family of divisors on $X$ $$f^i:Z^i\rightarrow C^i,$$ where $C^i$ is a curve with the following property: There are two points $P^i_0$ and $P^i_1$ on $C^i$ so that 
$$(f^i)^{-1}(P^i_0)=Z^i_0$$ and
$$(f^i)^{-1}(P^i_1)=Z^i_1,$$
and for any $P\in C^i$ $(f^i)^{-1}(P)$ is a divisor on $X$. It is possible to achieve this property in the following way. Assume $X$ is a projective variety. Let $X$ be the variety corresponding to the scheme $Proj(S)$ of finite type over $\overline{\Q}$, where $S=\oplus_{d=0}^{\infty}S_d$ is a graded ring generated by $S_1$ over $S_0=\overline{\Q}$. Let $I(Z^i_0)$ and $I(Z^i_1)$ be homogeneous ideals of homogeneous elements in $S$ that
vanish on the corresponding divisor. Let $s^i_0\in I(Z^i_0)$ and $s^i_1\in I(Z^i_1)$
be nonzero homogeneous elements of the same degree of the ideals so that $V((s^i_0))$ 
and $V((s^i_1))$ have no common component. Then for parameters $x$ and $y$ such that 
$(x:y)\in P^1(\overline{\Q})$, consider the divisors $V((xs^i_0 +ys^i_1))$. They form a family of divisors parametrized by $P^1(\overline{\Q})$. The family of divisors have the following property: above $(1:0)$ lies $V((s^i_0))\supset Z^i_0$ and 
above $(0:1)$ lies $V((s^i_1))\supset Z^i_1$. Thus possibly we need to add more components to our divisors.

Let $g^i:Z^i\rightarrow X$ be the map that sends  $(f^i)^{-1}(P)$ for $P\in C^i$
to the corresponding divisor on $X$. Let 
$$Z=Z^1\times_X\dots\times_X Z^n.$$ Let $p_i$ be the composition of 
$Z\rightarrow Z^i$ and $f^i:Z^i\rightarrow C^i$. 

Let us examine the fiber of 
$$g:Z\rightarrow X.$$ Note that $Z$ and $X$ are of the same dimension. If we assume that $X$ is projective then $g$ is proper and onto. Therefore there is a maximal Zariski open set
in $U\in X$ such that $g^{-1}(x)$ for $x\in U$ has finitely many elements.
We can restrict $U$ to a maximal Zariski open set $V$ so that 
$g:g^{-1}(V) \rightarrow V$ is \'etale. Let $Y=X-V$. Then $Y$ is a union of varieties of codimension at least $1$.

Let $d$ be the degree of the \'etale map. 

Let $$M_Z=g^{-1}(M\cap V).$$ Let $\omega\in\Ga(\Omega^n_X,X)$. Recall that
$X$ is of dimension $n$. Then
$$\int_{M_Z}\frac{1}{d}g^*\omega=\int_M \omega.$$ 
The boundary of $M_Z$ lies on divisors $D^i_0=g^{-1}(Z^i_0)$ and $D^i_1=g^{-1}(Z^i_1)$
inside $Z$. So the above integral is a period of a framed mixed motive.

We are going to show that the iterated integrals of $n$-forms from $\Ga(\Omega^n_X,X)$ on $n$-dimensional variety $X$ also has the structure of a period of a mixed motive.
For simplicity we will consider iteration of $2$-forms $\omega_1$ and $\omega_2$ on  a $2$-dimensional variety $X$. Let $M$ be admissible $2$-dimensional membrane on $X$,
where $M$ is the image of $\theta:[0,1]^2\rightarrow X(\C)$. Let $\a_i=\theta^*\omega_i$
for $i=1,2$. Consider the integral
$$\int_{0\leq x_1\leq x_2\leq 1; 0\leq y_2\leq y_1\leq 1}
\a_1(x_1,y_1)\wedge\a_2(x_2,y_2).$$
One can think of this integral as an integral defined on $X\times X$. 
We define $Z$ as above. Let $p_1:Z\rightarrow C^1$ and $p_2:Z\rightarrow C^2$ as above.
Pull back to $Z$ and consider this integral over $Z$ where we iterate $$\frac{1}{d}g^*\omega_i$$ for
$i=1,2$. The boundary for the domain of integration on $Z$ lies on certain divisors.
The divisors correspond to making one  inequality into an equality. The equality $0=x_1$ corresponds to $D^1_0\times Z$. The equality $x_2=1$ corresponds to 
$Z\times D^1_1$. The equality $0=y_2$ corresponds to $D^2_0\times Z$.
The equality $y_1=1$ corresponds to $Z\times D^2_1$. The equality $x_1=x_2$ corresponds
to $Z\times_{C^1}Z$ in $Z\times Z$. And the equality $y_2=y_1$ corresponds
to $Z\times_{C^2}Z$ in $Z\times Z$. Let $D_Z$ be the union of $D^1_0\times Z$, 
$Z\times D^1_1$,  $D^2_0\times Z$, $Z\times D^2_1$, $Z\times_{C^1}Z$ and 
$Z\times_{C^2}Z$ in $Z\times Z$.

Let $$M_2=\{(\theta(x_1,y_1),\theta(x_2,y_2))|
0\leq x_1\leq x_2\leq 1; 0\leq y_2\leq y_1\leq 1\}$$.
We are ready to show that $M_2$ is an admissible membrane on $X(\C)\times X(\C)$. We have that the boundary of $(g\times g)^{-1}{M_2}$ lies on union of varieties of codimension $1$ in $Z(\C)\times Z(\C)$. We have that $g$ is proper. So it is closed. Let $D'_X=(g\times g)(D_Z)$. Then the boundary of $M_2$ lies on $D$. Let $D_X(\C)$ be the Zariski closure in $X\times X$ of the boundary of $M_2$. Then $D_X(\C)\subset D'_X(\C)$.
All components of $D_X(\C)$ are components of $D'_X(\C)$. So $D_X(\C)$ is defined over $\overline{\Q}$. Denote it by $D_X$. The above iterated integral is represented by the motive $H_2(X\times X, D_X)$. It can be represented as framed mixed Hodge structure. The framing is given by non-trivial maps from $\bar{\Q}(0)$ to each copy of $\oplus \bar{\Q}(0)$ with one copy for each vertex of $M_2$, and the other map of the framing is given by $p_1^*\omega_1\wedge p_2^*\omega_2$ where 
$p_1$ and $p_2$ are the two projections of $X\times X$ to $X$.

Similarly, we can associate a framed motive to any iterated integral over an admissible membrane.

Note that if we vary an admissible membrane $M$ by a family of admissible
$M^t$ such that the boundary of $M^t$ stays on the same divisor of $X$ then the boundary $M^t_2$ will lie on the divisors of the boundary of $M_2$. Therefore the framed motive that correspond to an iterated integral over $M$ represent an integration over an element the fundamental $2$-groupoid of $X(\C)$ corresponding to an admissible membrane.
\section{Non-commutative 2 dimensional modular symbol for $SL_2$ over real quadratic rings}

Let $K$ be a real quadratic field and let ${\cal{O}}_K$ be its ring of integers.
Let $PGL_2({\cal{O}}_K)$ be the group of $2$ by $2$ invertible matrices with 
coefficients in ${\cal{O}}_K$ modulo units times the identity matrix. We shall define a non-abelian modular symbol
for a finite index subgroup $\Ga$ of $PGL_2({\cal{O}}_K)$. The space on which this group acts is a product of two copies of the upper half plane ${\cal{H}}\times{\cal{H}}={\cal{H}}_K$, where the action is determined by the two real embeddings of $K$ into $\R$. Identify each copy of ${\cal{H}}$ with the $2$ by $2$ positive definite quadratic forms, considered as matrices, modulo a positive real number times the identity matrix. We add point at infinity to ${\cal{H}}_K$ of the following type. Let 
$v^T=(\a,\b)$ be a non-zero vector over $K$. Let $P=v^Tv$. Consider the $2$ by $2$ positive semi-definite matrices of rank at least $1$ so that if a matrix is not positive definite then it is of the type $v^Tv$. Denote by
$\overline{{\cal{H}}}_K$ the space ${\cal{H}}_K$ with the points $P$ added at infinity.

The points $0=(0,0)$ $1=(1,1)$ and $\infty=(\infty,\infty)$ lie on $\overline{{\cal{H}}}$ which is a complex space and can be embedded diagonally
into $\overline{{\cal{H}}}_K$. Choose a point $P\in P^1(K)$ on the boundary of $\overline{{\cal{H}}}_K$. Let $g_0$, $g_1$ and $g_2$ be elements of $\Ga$. Consider the points $P_0=g_0P$, $P_1=g_1P$, $P_2=g_2P$. There exists unique $g\in PGL_2(K)$ such that $gP_0=0$, $gP_1=1$ and $gP_2=\infty$. So $P_0$, $P_1$ and $P_2$ lie on 
$g^{-1}\overline{{\cal{H}}}$ which is again a space with complex structure. In the quotient $\Ga\backslash\overline{{\cal{H}}}_K$ $g^{-1}\overline{{\cal{H}}}$
becomes an algebraic curve. Consider the geodesic that connects $P_0$ and $P_1$ that lies on $g^{-1}\overline{{\cal{H}}}$. Call it $\ga_{P_0P_2P_1}$. Let $g_3\in \Ga$. And let $P_3=g_3P$. Take the unique $g'\in PGL_2(K)$ such that $g'P_0=0$, $g'P_1=1$ and  
$g'P_3=\infty$. Let $\ga_{P_0P_3P_1}$ be the geodesic that connects $P_0$ and $P_1$ so
that it lies on $g'^{-1}\overline{{\cal{H}}}$.

We shall iterate holomorphic modular 2-forms which are modular with respect to $\Ga$. 
The membrane of integration will have boundaries $\ga_{P_0P_2P_1}$ and $\ga_{P_0P_3P_1}$. The generating series for the iterated integral are homotopy invariant if we keep the boundary fixed. We can vary the membrane of integration even more. We can vary $\ga_{P_0P_2P_1}$ so that the end points stay fixed and the variation lies on $g^{-1}\overline{{\cal{H}}}$. Then the generating series will not change, since integration of holomorphic two forms on a 1 dimensional complex space $g^{-1}\overline{{\cal{H}}}$ gives zero. When we pass to the quotient 
$\Ga\backslash\overline{{\cal{H}}}_K$ we can vary the geodesics $\ga_{P_0P_2P_1}$ and $\ga_{P_0P_3P_1}$ so that the variation of $\ga_{P_0P_2P_1}$ lies in the algebraic curve which is image of $g^{-1}\overline{{\cal{H}}}$ and the geodesic $\ga_{P_0P_3P_1}$ lies 
on the algebraic curve which is image of $g'^{-1}\overline{{\cal{H}}}$.  In this way we can iterate cusp forms for $\Ga$ which come from $H^0(X(\Ga),V)$, where $X(\Ga)$ is the corresponding Hilbert modular surface defined over $\overline{\Q}$. Using the constructions in the section on motivic interpretation of iterated integrals of cusp forms, we obtain period of framed motives. Note that the membrane that we take here is an admissible membrane as defined in the section on motivic interpretation of iterated integrals. In this way we can iterate cusp forms for $\Ga$. And we obtain non-abelian modular symbol for $\Ga$ defined as the generating series $$J_M(\omega_1,\dots,\omega_n),$$
where $M$ is the membrane defined above, via geodesics, and $\omega_i$ are cusp forms.

We would like to formulate one conjecture on the coproduct of the framed motives corresponding to iterated integrals. Let $M$ be a membrane of the above type on a Hilbert modular surface. Let 
$$I=I_M(\omega_1,\dots,\omega_n,\si_1,\si_2)$$
be an iterated integral of cusp forms $\omega_i$ for $i=1,\dots,n$ with permutations
$\si_1$ and $\si_2$ of $\{1,2\dots,n\}$. For each iterated integral $I$, let $\tilde{I}$ be the corresponding framed motive. Let
$$\Delta(I)=\oplus_i I'_i\otimes I''_i$$ be the coproduct of the iterated integral $I$, as defined in section $1$. Let $\tilde{\Delta}$ be the coproduct of framed motives.
\bcon{6.2}
With the above notation 
$$\tilde{\Delta}\tilde{I}=\oplus_i \tilde{I'_i}\otimes \tilde{I''_i},$$
where $\Delta(I)=\oplus_i I'_i\otimes I''_i$ is the coproduct of iterated integrals.
Therefore the map from iterated integrals of the above type to framed motives is a morphism of Hopf algebras.
\econ
We are going to examine multiple $L$-functions associated to Hilbert cusp forms. Let
$\Ga$ be a finite index subgroup of $GL_2({\cal{O}}_K)$ or of $SL_2({\cal{O}}_K)$.
Let $p$ be the projection from $\overline{\cal{H}}\times \overline{\cal{H}}$ to the Hilbert modular surface $X(\Ga)$. Let $M$ be the geodesic membrane
in $\overline{\cal{H}}\times \overline{\cal{H}}$ consisting of geodesics connecting $0$ and $\infty$, and bounded by $\ga(0,1,\infty)$ and $\ga(0,u,\infty)$ where $u$ is a generator of the group of units in $K$. The membrane $p(M)$ is an admissible membrane because its boundaries lie on algebraic curves in $X(\Ga)$. We are going to define the multiple $L$-functions as iterated Mellin transforms, where the integration is taken over $M$. Let 
$$L(\omega_1,\dots,\omega_n,s_1,\dots,s_n,\si_1,\si_2)=
I_M(\omega_1(-z_1z_2)^{s_1-1},\dots,
\omega_n(-z_1z_2)^{s_1-1},\si_1,\si_2).$$
Note that $$I_M(\omega_1,dz_1\wedge dz_2,dz_1\wedge dz_2,
\omega_2,(1),(1))=L(\omega_1,\omega_2,3,1,(1),(1)),$$
where $(1)$ it the identity permutation. Therefore the positive integral values of the multiple $L$-functions are periods of framed motives.
\section{Multiple completed Dedekind zeta functions}

\subsection{Multiple completed Rieman zeta function}

We give an application of this process to the iteration of the completed Riemann zeta function. For this definition we use iterated integrals over a path as defined by Chen \cite{Ch}. Recall the Riemann zeta function is
$$\zeta(s)=\sum_{n=1}^\infty \frac{1}{n^s}.$$
Consider the following Jacobi theta function
$$\theta(z)=
\sum_{n\in\Z}e^{\pi i(n^2z)}.$$
Certain Mellin transform of $\theta(z)-1$ gives us $\zeta(s)$ times a gamma factor.
The Mellin transform that we will consider here is 
$$\hat{\zeta}(s)=
\int_\ga (\theta(z)-1)(-iz)^{\frac{s}{2}}\frac{dz}{z}.$$
where $\ga$ is the geodesic connecting $0$ and $\infty$ in the upper half plane. Define
$$\Ga'(2)=\left\{
A\in SL_2({\cal{O}}_K)|A \equiv 
\left[\begin{tabular}{cc}
$\pm 1$ & $0$\\
$0$ & $\pm 1$
\end{tabular}\right]
\mbox{ or }
\left[\begin{tabular}{cc}
$0$ & $\pm 1$\\
$\mp 1$& $0$
\end{tabular}\right]
\mbox{ mod } 2
\right\}.$$
Let $Y={\cal{H}}/\Ga'(2)$. Then $\theta(z)$ can be defined as a section of a line bundle $\cal{L}$ of $Z$ where $Z\rightarrow Y$ is a finite morphism.
We are going to use that $\theta(z)$ is defined algebraically on $Z$ over $\Q$. If we believe that, we are going to show that iterated theta function gives period over framed motives which we are going to call multiple Dedekind zeta functions.

We can consider the iterated integrals of $(\theta(z)-1)dz$ and $dz$ as period of mixed Hodge structures which are associated to algebraic varieties defined over $\Q$. 

\bde{6.9} Define multiple completed zeta function to be
$$\hat{\zeta}(s_1,\dots,s_d,\si_1,\si_2)=
I_{\ga}((\theta(z)-1)(-iz)^{s_1},\dots,
(\theta(z)-1)(-iz)^{s_d}).$$
\ede

>From the above comments we have that $$\hat{\zeta}(n_1,\dots,n_d)$$ has motivic interpretation when $n_i$ are positive even integers. For example 
$$\hat{\zeta}_K(4,2,(1))=
I_{\ga}((\theta(z)-1)dz,dz,
(\theta(z)-1)dz).$$
\subsection{Multiple completed Dedekind zeta functions for imaginary quadratic field of class number $1$}

We give an application of this process to the iteration of Dedekind zeta function for imaginary quadratic field $K$ of class number one. Let $K$ be a imaginary quadratic field of class number one. And let ${\cal{O}}_K$ be the ring of integers in $K$. Then 
$$\zeta_K(s)=\sum_{(\a)\subset {\cal{O}}_K,(\a)\neq (0),} \frac{1}{N((\a))^s}.$$

We define a theta function whose Mellin transform gives $\zeta_K(s)$. The group
$SL_2({\cal{O}}_K)$ acts on the upper half space. We can identify the upper half space with the space $2$ by $2$ positive definite Hermitian matrices. Let $x$ be a positive definite Hermitian matrix. Consider the matrix $x+\bar{x}$. It is a $2$ by $2$ a positive definite symmetric matrix. We identify the space of $2$ by $2$ positive definite symmetric matrices with the upper half plane sitting inside $\C$. Denote by $[x]$ the complex number in the upper half plane that corresponds to $x+\bar{x}$.
Define
$$\theta_K(x)=
\sum_{\a\in{\cal{O}}_K}e^{\pi i(N(\a)[x])}.$$
Certain Mellin transform of $\theta_K(x)-1$ gives us $\zeta_K(s)$ times a gamma factor.
The Mellin transform that we will consider here is 
$$\hat{\zeta}_K(s)=\int_{\ga} (\theta_K(x)-1)(-i[x])^s\frac{d[x]}{[x]}.$$
where $ga$ is the geodesic connecting $0$ and $\infty$ in the upper half plane.

The modular form $\theta_K(x)$ is modular with respect to a group $\Ga'_K$, where
$$\Ga'_K=\left\{
A\in SL_2({\cal{O}}_K)|A \equiv 
\left[\begin{tabular}{cc}
$u$ & $0$\\
$0$ & $u^{-1}$
\end{tabular}\right]
\mbox{ or }
\left[\begin{tabular}{cc}
$0$ & $u$\\
$-u^{-1}$& $0$
\end{tabular}\right]
\mbox{ mod } Z[\sqrt{-D}]
\right\},$$
where $u$ varies through all invertible elements and $-D$ is the discriminant of the field $K$.

Consider the geodesic plane in the upper half space that passes through $0$, $1$ and $\infty$. We can identify this geodesic plane with the upper half plane. It correspond to taking the real part of the positive definite $2$ by $2$ Hermitian matrices. 
Restrict the function $\theta_K(x)$ to the upper half plane that passes through 
$0$, $1$ and $\infty$. Denote this function by $\theta_K(z)$, where $z$ varies in the upper half plane. Then $\theta_K(z)$ is a modular function with respect to the group 
$$\Ga'_K\cap SL_2(\Z)=SL_2(\Z).$$

Let $Y={\cal{H}}/SL_2(\Z)$. Then $\theta_K(z)$ can be defined as a section of a line bundle $\cal{L}$ of $Z$ where $Z\rightarrow Y$ is a finite morphism.
We are going to use that $\theta_K(z)$ is defined algebraically on $Z$ over $\bar{\Q}$. If we believe that, we are going to show that iterated theta function gives period of motives which we are going to call multiple Dedekind zeta functions.

\bde{6.10} Define multiple completed Dedekind zeta function to be
$$\hat{\zeta}_K(s_1,\dots,s_d,\si_1,\si_2)=
I_{\ga}((\theta_K(z)-1)(-iz)^{s_1},\dots,
(\theta_K(z)-1)(-iz)^{s_d}).$$
\ede

>From the above coments we have that $$\hat{\zeta}_K(n_1,\dots,n_d)$$ has motivic interpretation when $n_i$ are positive integers. For example 
$$\hat{\zeta}_K(2,1)=
I_{\ga}((\theta_K(z)-1)dz,dz,(\theta_K(z)-1)dz).$$

\subsection{Multiple completed Dedekind zeta functions for real quadratic field of class number $1$}

Another application of this process is the iteration of Dedekind zeta function for real quadratic field $K$ of class number one. Let $K$ be a real quadratic field of class number one. And let ${\cal{O}}_K$ be the ring of integers in $K$. Then 
$$\zeta_K(s)=\sum_{(\a)\subset {\cal{O}}_K,(\a)\neq (0),} \frac{1}{N((\a))^s}.$$
Let $\si_1$ and $\si_2$ be the two real embeddings of $K$ into $\R$. Given $\a\in K$ let $\a_1=\si_1(\a)$ and $\a_2=\si_2(\a)$. We define the following function
$$\theta_K(z_1,z_2)=
\sum_{\a\in{\cal{O}}_K}e^{\pi i(\a_1^2z_1 + \a_2^2z_2)}.$$
Certain Mellin transform of $\theta_K(z_1,z_2)-1$ gives us $\zeta_K(s)$ times two gamma factors.
The Mellin transform that we will consider here is 
$$\hat{\zeta}_K(s,t)=\int\int_M \frac{1}{2}(\theta_K(z_1,z_2)-1)(-iz_1)^{\frac{s}{2}}(-iz_2)^{\frac{t}{2}}
\frac{dz_1}{z_1}\wedge \frac{dz_2}{z_2}.$$
The membrane $M$, over which we integrate, consists of family of geodesics connecting $0$ and $\infty$. The boundary of $M$ are two geodesics. One is $\ga_{0,1,\infty}$ and the other $\ga_{0,u^2,\infty}$, where $u$ is a generator for the units in $K$. Recall 
$\ga_{0,1,\infty}$ is the geodesic that lies on the $\Delta(\cal{H})\subset \cal{H}\times\cal{H}$ note that $0,1,\infty$ are boundary point of $\Delta(\cal{H})$. 
And $\ga_{0,u,\infty}$ is the geodesic that lies on $A\Delta(\cal{H})$, where
$A$ is the matrix in $GL_2({\cal{O}}_K)$ given by
$$A_{u^2}=\left[\begin{tabular}{cc}
$u$ & $0$\\
$0$ & $u^{-1}$
\end{tabular}\right].$$
The reason for this choice of $M$ is the following. In order to obtain the Dedekind zeta function, we need to integrate over the plane consisting of all geodesics between $0$ and $\infty$ modulo the action of the units. The action of $u$ is given by
$$A_u=\left[\begin{tabular}{cc}
$u$ & $0$\\
$0$ & $1$
\end{tabular}\right],$$ 
where we consider $A_u$ as an element of $PGL_2({\cal{O}}_K)$.
We have the fraction $1/2$ in front of the theta function because we factor over geodesics between $0$ and $\infty$ by the group generated by $u^2$, where $u$ is a generator of the group of units. We have defined $M$ in a such a way. The reason for quotienting by $u^2$ is that $A_{u^2}$, as defined above, is an element of the group 
$\Ga'_K(2)$ with respect to which $\theta_K(z_1,z_2)$ is a modular function. We define 
$\Ga'_K(2)$ below.
Note that $$\hat{\zeta}_K(s,s)=\pi^{-s}\Ga\left(\frac{s}{2}\right)^2\zeta_K(s).$$

The modular form $\theta(z_1,z_2)$ is modular with respect to $\Ga'_K(2)$ where
$$\Ga'_K(2)=\left\{
A\in SL_2({\cal{O}}_K)|A \equiv 
\left[\begin{tabular}{cc}
$u$ & $0$\\
$0$ & $u^{-1}$
\end{tabular}\right]
\mbox{ or }
\left[\begin{tabular}{cc}
$0$ & $u$\\
$-u^{-1}$& $0$
\end{tabular}\right]
\mbox{ mod } 2
\right\},$$
where $u$ varies through all invertible elements modulo $2$.

Let $Y={\cal{H}}\times{\cal{H}}/\Ga'_K(2)$. Then $\theta_K(z_1,z_2)$ can be defined as a section of a line bundle $\cal{L}$ of $Z$ where $Z\rightarrow Y$ is a finite morphism.
We are going to use that $\theta_K(z_1,z_2)$ is defined algebraically on $Z$ over $\bar{\Q}$. If we believe that, we are going to show that iterated theta function gives period over framed motives which we are going to call multiple Dedekind zeta functions.

By the construction in the beginning of this section we know that the image of $M$ under the projection to $Y$ is an admissible membrane. That is, its boundaries lie on algebraic curves in $Y$ defined over $\bar{\Q}$. Therefore we can consider the iterated integrals of $1/2(\theta(z_1,z_2)-1)dz_1\wedge dz_2$ and $dz_1\wedge dz_2$ as period of mixed Hodge structures which are associated to algebraic verieties defined over $\bar{\Q}$. 

\bde{6.11} Define multiple completed Dedekind zeta function to be
$$\hat{\zeta}_K(s_1,\dots,s_d,\si_1,\si_2)=
I_M(\frac{1}{2}(\theta(z_1,z_2)-1)(-z_1z_2)^{s_1},\dots,
\frac{1}{2}(\theta(z_1,z_2)-1)(-z_1z_2)^{s_d},\si_1,\si_2).$$
\ede

>From the above comments we have that $$\hat{\zeta}_K(n_1,\dots,n_d,\si_1,\si_2)$$ has motic interpretation when $n_i$ are positive even integers, for any permutations $\si_1$ and $\si_2$. For example 
$$\hat{\zeta}_K(4,2,(1),(12))=
I_M(\frac{1}{2}(\theta(z_1,z_2)-1)dz_1\wedge dz_2,dz_1\wedge dz_2,
\frac{1}{2}(\theta(z_1,z_2)-1)dz_1\wedge dz_2,(1),(13)).$$
\subsection{Multiple completed Dedekind zeta functions for CM fields of degree $4$ over $\Q$ and of class number $1$}

We give an application of this process to the iteration of Dedekind zeta function for a CM field $K$ of degree $4$ over $\Q$ and of class number one. Let $K$ be such a field. And let ${\cal{O}}_K$ be the ring of integers in $K$. Then 
$$\zeta_K(s)=\sum_{(\a)\subset {\cal{O}}_K,(\a)\neq (0),} \frac{1}{N((\a))^s}.$$

We define a theta function whose Mellin transform gives $\zeta_K(s)$. The group
$SL_2({\cal{O}}_K)$ acts on a product of two upper half spaces. We can identify the upper half space with the space $2$ by $2$ positive definite Hermitian matrices. Let $x$ be a positive definite Hermitian matrix. Consider the matrix $x+\bar{x}$. It is a $2$ by $2$ a positive definite symmetric matrix. We identify the space of $2$ by $2$ positive definite symmetric matrices with the upper half plane sitting inside $\C$. Denote by $[x]$ the complex number in the upper half plane that corresponds to $x+\bar{x}$. Given an element $\a$ in $K$, let $\si_1$ and $\si_2$ be two non-conjugate complex embeddings of $K$ into $\C$. Let also $\a_i=\si_i(\a)$ for $i=1,2$. Let also $(x_1,x_2)$ be a point in the product of two upper half spaces. Define
$$\theta_K(x_1,x_2)=
\sum_{\a\in{\cal{O}}_K}e^{\pi i(\a_1\bar{\a_1}[x_1]+\a_2\bar{\a_2}[x_2])},$$
where $\bar{\a_i}$ is the complex conjugate of $\a_i$. 
Certain Mellin transform of $\theta_K(x_1,x_2)-1$ gives us $\zeta_K(s)$ times two gamma factors.
The Mellin transform that we will consider here is 
$$\hat{\zeta}_K(s,t)=\int\int_{M} \frac{1}{d}(\theta_K(x_1,x_2)-1)(-i[x_1])^s(-i[x_2])^t\frac{d[x_1]}{[x_1]}\wedge \frac{d[x_2]}{[x_2]},$$
The choice of $M$ is the following. In order to obtain the Dedekind zeta function we need to integrate over the plane consisting of all geodesics between $0$ and $\infty$ modulo the action of the units. The action of $u$ is given by
$$A_u=\left[\begin{tabular}{cc}
$u$ & $0$\\
$0$ & $1$
\end{tabular}\right],$$ 
where we consider $A_u$ as an element of $PGL_2({\cal{O}}_K)$.
The modular form $\theta_K(x_1,x_2)$ is modular with respect to a group $\Ga'_K$, where
$$\Ga'_K=\left\{
A\in PSL_2({\cal{O}}_K)|A \equiv 
\left[\begin{tabular}{cc}
$u$ & $0$\\
$0$ & $u^{-1}$
\end{tabular}\right]
\mbox{ or }
\left[\begin{tabular}{cc}
$0$ & $u$\\
$-u^{-1}$& $0$
\end{tabular}\right]
\mbox{ mod } R
\right\},$$
where $$R=\{\a\in{\cal{O}}_K|\a_i+\bar{\a_i}\in(2),\mbox{ for } i=1,2\},$$ and $u$ is a unit in ${\cal{O}}_K$. Let $d_0$ be the smallest positive integer such that $A_{u_0^d}\in\Ga'_K$ and $u_0$ is the generator of the unit group of the totally real subfield. Let $U_K$ the group of units in $K$ and let $<u^d>$ be the subgroup generated by $u^d$. Define $$d=\#|U_K/\mu(K)<u^d>|.$$
Take $M$ to be the subset of all geodesics between $0$ and $\infty$ which is bounded between $\ga_{0,1,\infty}$ and $\ga_{0,u_0^{d_0},\infty}$. Recall $\ga_{0,1,\infty}$ is the geodesic that connects $0$ and $\infty$ which lies in the geodesic plane that passes through $0$, $1$ and $\infty$. (This geodesic plane is an upper half plane.)

Consider the geodesic plane in the upper half space that passes through $0$, $1$ and $\infty$. We can identify this geodesic plane with the upper half plane. It corresponds to taking the real part of the positive definite $2$ by $2$ Hermitian matrices. 
Restrict the function $\theta_K(x_1,x_2)$ to the upper half plane that passes through 
$0$, $1$ and $\infty$. Denote this function by $\theta_K(z_1,z_2)$, where $z$ varies in the upper half plane. Then $\theta_K(z)$ is a modular function with respect to the group 
$$\Ga'_K\cap PSL_2(\Z)=\Ga''_K.$$

Let $Y={\cal{H}}\times {\cal{H}}/\Ga''_K$. Then $\theta_K(z_1,z_2)$ can be defined as a section of a line bundle $\cal{L}$ of $Z$ where $Z\rightarrow Y$ is a finite morphism.
We are going to use that $\theta_K(z_1,z_2)$ is defined algebraically on $Z$ over $\bar{\Q}$. If we believe that, we are going to show that iterated theta function gives period of motives which we are going to call multiple Dedekind zeta functions.
Note that $M$ is an admissible membrane. It can be projected to the algebraic variety
$Y$ and the boundaries of $M$ lie on images of upper half planes, which map under the projection modulo  $\Ga''_K$ to algebraic curves.
\bde{6.12} Define multiple completed Dedekind zeta function to be
$$\hat{\zeta}_K(s_1,\dots,s_d,\si_1,\si_2)=
I_{M}(\frac{1}{d}(\theta_K(z_1,z_2)-1)(-z_1z_2)^{s_1},\dots,
\frac{1}{d}(\theta_K(z_1,z_2)-1)(-z_1z_2)^{s_d}).$$
\ede

>From the above coments we have that $$\hat{\zeta}_K(n_1,\dots,n_d)$$ has motivic interpretation when $n_i$ are positive integers. For example 
$$\hat{\zeta}_K(2,1,(1),(12))=
I_{M}(\frac{1}{d}(\theta_K(z_1,z_2)-1)dz_1\wedge dz_2,dz_1\wedge dz_2,
\frac{1}{d}(\theta_K(z_1,z_2)-1)dz_1\wedge dz_2,(1),(13)).$$
\section{Final remarks}
The iteration over membranes that we consider can easily be generalized to higher dimensions. Then one can construct non-abelian symbols for $GL_m(\Z)$. However, there is a difficult problm in combinatorics that needs to be solves, in order to show that the generating series of the interated integrals over a membrane satisfy a cocycle condition. But the actual generating series can be defined easily.

About the iterated completed Dedekind zeta functions, I would like to mension that for $K=\Q$ and for $K$ inaginary quadraticfield of class number one, we obtain mixed Tate motive. I would state a conjecture. 
\bcon{6.15}
For $K$ a totally real field of class number one or a CM field of degree $4$ and of clas number one, the motives corresponding to the iterated completed Dedekind zeta fuctions are mixed Tate motives.
\econ

\renewcommand{\em}{\textrm}

\begin{small}

\renewcommand{\refname}{ {\flushleft\normalsize\bf{References}} }
    
\end{small}

\end{document}